\documentclass[twoside,11pt]{article}
\title{} \author{} \date{}
\usepackage{latexsym,amssymb,times}
\input amssym.def
\newtheorem{te}{Theorem}[section]
\newtheorem{prop}[te]{Proposition}

\newtheorem{fac}[te]{Fact}
\newtheorem{lem}[te]{Lemma}

\newtheorem{rem}[te]{Remark}
\newtheorem{ex}[te]{Example}
\newtheorem{que}[te]{Question}

\def\dok{\noindent{\bf Proof. }}
\def\kdok{\hfill $\Box$ \par \vspace*{2mm} }
\def\a{\alpha}
\def\b{\beta}
\def\g{\gamma}
\def\d{\delta}
\def\f{\varphi}

\def\o{\omega}
\def\k{\kappa}
\def\l{\lambda}
\def\r{\rho}
\def\s{\sigma}

\def\h{{\mathfrak h}}
\def\c{{\mathfrak c}}

\def\S{{\mathbb S}}
\def\P{{\mathbb P}}
\def\Q{{\mathbb Q}}
\def\B{{\mathbb B}}
\def\N{{\mathbb N}}
\def\X{{\mathbb X}}
\def\Y{{\mathbb Y}}

\def\A{{\mathbb A}}

\def\M{{\mathbb M}}

\def\BR{{\mathbb R}}
\def\BI{{\mathbb I}}
\def\BF{{\mathbb F}}
\def\BL{{\mathbb L}}
\def\BK{{\mathbb K}}
\def\BH{{\mathbb H}}
\def\BE{{\mathbb E}}

\def\CP{{\mathcal P}}

\def\CD{{\mathcal D}}

\def\I{{\mathcal I}}

\def\L{{\mathcal L}}

\def\la{\langle}
\def\ra{\rangle}

\def\dom{\mathop{\mathrm{dom}}\nolimits}

\def\id{\mathop{\mathrm{id}}\nolimits}
\def\Borel{\mathop{\mbox{Borel}}\nolimits}

\def\Scatt{\mathop{\mbox{\rm Scatt}}\nolimits}

\def\Emb{\mathop{\rm Emb}\nolimits}

\def\Aut{\mathop{\rm Aut}\nolimits}

\def\Fin{\mathop{\rm Fin}\nolimits}
\def\Coll{\mathop{\rm Coll}\nolimits}

\def\sm{\mathop{\rm sm}\nolimits}
\def\sq{\mathop{\rm sq}\nolimits}
\def\ro{\mathop{\rm ro}\nolimits}

\def\Lim{\mathop{\mbox{Lim}}\nolimits}
\def\rp{\mathop{\rm rp}\nolimits}

\def\Borel{\mathop{\rm Borel}\nolimits}
\def\Mod{\mathop{\rm Mod}\nolimits}

\def\ar{\mathop{\rm ar}\nolimits}
\def\cf{\mathop{\rm cf}\nolimits}

\def\Col{\mathop{\rm Col}\nolimits}

\def\Pa{\mathop{\rm Pa}\nolimits}

\def\Sym{\mathop{\rm Sym}\nolimits}

\def\Sym{\mathop{\rm Sym}\nolimits}

\def\otp{\mathop{\rm otp}\nolimits}
\def\Copy-dense{\mathop{\mathrm{Copy-dense}}\nolimits}
\def\LO{\mathop{\rm LO}\nolimits}

\def\cc{\mathop{\rm cc}\nolimits}
\def\sh{\mathop{\rm sh}\nolimits}
\def\ord{\mathop{\rm Ord}\nolimits}
\def\RO{\mathop{\rm RO}\nolimits}

\begin{document}
\thispagestyle{plain}
\begin{center}
           {\large \bf {\uppercase{Copies of monomorphic structures}}}
\end{center}
\begin{center}
{\bf Milo\v s S.\ Kurili\'c}\footnote{Department of Mathematics and Informatics, Faculty of Sciences, University of Novi Sad, \\
         Trg Dositeja Obradovi\'ca 4, 21000 Novi Sad, Serbia. e-mail: milos@dmi.uns.ac.rs}
\end{center}
\begin{abstract}
\noindent
The poset of copies of a relational structure ${\mathbb X}$ is the partial order $\langle {\mathbb P} ({\mathbb X}) ,\subset \rangle$,
where ${\mathbb P} ({\mathbb X})=\{ Y\subset X: {\mathbb Y} \cong {\mathbb X}\}$.
Investigating the classification of structures
related to isomorphism of the Boolean completions ${\mathbb B}_{\mathbb X} ={\mathop{\rm ro}\nolimits}({\mathop{\rm sq}\nolimits} ({\mathbb P} ({\mathbb X}) ))$
we extend the results concerning linear orders to the class of
structures definable in linear orders by first-order $\Sigma _0$-formulas (monomorphic structures).
So, ${\mathbb B}_{\mathbb X} \cong {\mathbb B}_{\mathbb L}$ holds for some linear order ${\mathbb L}$,
if ${\mathbb X}$ is definable in a $\sigma$-scattered (in particular, countable) or additively indecomposable linear order.
For example, ${\mathbb B}_{\mathbb X} \cong {\mathop{\rm ro}\nolimits}({\mathbb S} )$, where ${\mathbb S}$ is the Sacks forcing,
whenever ${\mathbb X}$ is a non-constant structure chainable by a real order type containing a perfect set.
\end{abstract}
\section{Introduction}\label{S1}
If $\X$ is a relational structure
and $\Emb (\X )$ is the set of its self-embeddings,
then the partial order $\la \P (\X ) ,\subset \ra$,
where $\P (\X ):=\{ f[X]: f\in \Emb (\X )\}$,
is called the {\it poset of copies of} $\X$
and $\sq (\P (\X))$ and $\B_\X :=\ro(\sq (\P (\X) ))$ denote its separative quotient and Boolean completion.
The implications
\begin{equation}\label{EQ004}
\X \cong \Y \Rightarrow \P(\X )\cong \P(\Y) \Rightarrow \sq \P(\X )\cong \sq \P(\Y) \Rightarrow \B_\X \cong \B_\Y
\end{equation}

\vspace*{-1mm}
\noindent
express the relationship between different ``similarity relations" in the class of relational structures
and each of them provides a classification of structures \cite{Ktow,KDif}.
The last similarity, $\B_\X \cong \B_\Y$, induces the coarsest classification of structures and
it is equivalent to the forcing equivalence of posets of copies, namely (see \cite{KDif})
\begin{equation}\label{EQ005}
\P(\X )\equiv_{forc} \P(\Y)\Leftrightarrow \B_\X \cong \B_\Y.
\end{equation}

\vspace*{-1mm}
\noindent
So, the corresponding classification of structures and their embedding monoids can be explored using forcing-theoretic methods and related facts from set theory.

The aim of this paper is to extend the known results
concerning linear orders to the class of relational structures
which are definable in linear orders by first-order formulas without quantifiers--the class of infinite monomorphic structures.
It seems that such extension achieves a natural bound,
since the class of monomorphic structures is closed under substructures, bi-embeddability,
$\Sigma _0$-definability and elementary equivalence \cite{KVau}.

Roughly, our strategy is the following.
If $\X$ is a monomorphic structure $\Sigma _0$-definable in a linear order $\BL$,
then (regarding our goal)
it is natural to ask whether the algebras $\B _\X$ and $\B _\BL$ are isomorphic.
Trivially, $\B _\X \not\cong\B _\BL$ is possible:
if $\Q =\la Q, <\ra$ is the rational line
and $\X =\la Q, \r\ra$, where $\r$ is the binary relation defined in $\Q$ by the formula $\f (u,v):=u\neq v$
(that is, $\X$ is the complete graph $\la Q, Q^2\setminus \Delta _Q \ra$),
then  $\B _\X \cong \ro (\P (\o )/\Fin)$, while $\B _\Q \cong \ro (\S \ast \pi)$, where $\S$ is the Sacks forcing (see Theorems \ref{T010}(b) and \ref{T011}(a)).
Fortunately, for a monomorphic structure $\X$
Fact \ref{T8090} describes the set $\L _\X$ of {\it all} linear orders $\BL=\la X,<\ra$ in which $\X$ is  $\Sigma _0$-definable;
thus, a better question is:
\begin{que}\label{Q1}
If $\X$ is a monomorphic structure,
is there a linear order $\BL\in \L _\X$ such that $\B _\X \cong\B _{\BL}$?
\end{que}
Then knowing $\B_\BL$ we obtain $\B _\X$ for free.

Section \ref{S2} contains necessary definitions and facts.
In Section \ref{S3} we consider the simplest monomorphic structures, called {\it constant} by  Fra\"{\i}ss\'{e},
and the structures from a larger class ({\it copy-dense structures}, in our terminology);
here, the main role is played by the algebras of the form $P(\k )/[\k]^{<\k}$, which are consistently isomorphic to collapsing algebras.
Question \ref{Q1} is considered in Section \ref{S4};
it turns out that we have the answer `Yes'
if the collection $\L _\X$ contains a $\hookrightarrow$-minimal element, or an additively indecomposable, or a $\s$-scattered linear order (by Laver's theorem).
Consequently, we have the answer `Yes' for the structures $\Sigma _0$-definable in well orders,
and regarded in Section \ref{S5}; here we detect a minimal element of the quasi-order $\la \L _\X ,\hookrightarrow\ra$ (existing by Laver's theorem)
and obtain a representation of $\B _\X$ in terms of the Cantor normal form for the related ordinal.
For countable monomorphic structures considered in Section \ref{S6} we have the answer `Yes' again;
from the previous results and the known facts from \cite{Kscatt,Kord,KurTod}
it follows that the main roles are played by the Sacks forcing and the algebras $P(\o ^\d)/\I _{\o ^\d}$, where $\d <\o _1$.
In particular, under CH for each countable chainable structure $\X$ we have
$$
\B _\X \cong
               \left\{\begin{array}{cl}
                 \ro (P(\o )/\Fin),                   & \mbox{if $\X$ is chained by a scattered l.\ o.,} , \\
                 \ro (\S \ast \dot{(P(\o )/\Fin)^+}), & \mbox{if $\X$ is chained by a non-scattered l.\ o.}
               \end{array}
               \right.
$$
In Section \ref{S7} we first visit the large zoo of structures definable in uncountable ordinals;
using the results from \cite{Kord1} we describe the algebras $\B _\X$ in several models of ZFC;
here, important roles are played by collapsing algebras.
In the sequel we consider the structures chained by uncountable real types
and confirm the fundamental role of the Sacks forcing in this context,
showing that (in ZFC)
$$
\B _\X \cong \ro (\S )\cong \ro(\Borel (\BR )/[\BR]^{\leq \o}),
$$
whenever $\X$ is a non-constant relational structure chainable by a real order type containing a perfect set.
\section{Preliminaries}\label{S2}
\paragraph{Monomorphic structures}
Let $L=\la R_i:i\in I\ra$ be a relational language, where $\ar (R_i)=n_i$, for $i\in I$, and let $\X =\la X,\la \r _i:i\in I\ra\ra =\la X,\bar{\r}\ra$ be an $L$-structure.
By $\Pa (\X)$ we denote the set of all partial automorphisms (i.e., isomorphisms between substructures) of $\X$,
by $\LO_ X$ the set of all linear orders on the set $X$
and by $\CD _{\Sigma _0}(\LO_X)$ the set of all $L$-structures which are definable in linear orders from $\LO _X$ by the first order $L_b$-formulas without quantifiers and parameters.
Thus, $\X\in \CD _{\Sigma _0}(\LO)$ iff there is a linear order $<$ on $X$ and for each $i\in I$ there is an $L_b$-formula $\f (v_0,\dots,v_{n_i-1})$ as above such that
$\r _i=\{ \bar x \in X^{n_i}: \la X,<\ra \models \f _i[\bar x]\}$.

The structure $\X$ is called {\it monomorphic} iff each two finite substructures of $\X$ of the same size are isomorphic;
$\X$ is  called {\it chainable} if there is a linear order $<$ on $X$ such that $\Pa (\la X,< \ra )\subset \Pa (\X)$.
The following statement was proved by Fra\"{\i}ss\'{e} (see \cite{Fra}) for finite languages and by Pouzet for arbitrary languages \cite{Pou}.
\begin{fac}\label{T8067}
An infinite relational structure $\X$ is monomorphic iff it is chainable iff $\,\X\in \CD _{\Sigma _0}(\LO _X)$.
\end{fac}
\begin{ex}\label{EX8008}\rm
If $\la X,<\ra$ is a linear order, then the {\it betweness relation}, $D_{\f _b}\subset X^3$,
is defined by the formula $\f_b:=(v_0<v_1<v_2) \lor (v_2<v_1<v_0)$;
the formula $\f_c:=(v_0<v_1<v_2) \lor (v_1<v_2<v_0)\lor (v_2<v_0<v_1)$ defines the {\it cyclic relation}, $D_{\f _c}\subset X^3$;
and the {\it separation relation}, $D_{\f _s}\subset X^4$, is defined by a formula $\f _s (v_0,v_1,v_2,v_3)$
saying that $v_0,v_1,v_2$ and $v_3$ are different and the pair $\{ v_0,v_2 \}$ separates the pair $\{ v_1,v_3 \}$ (see \cite{Fra}).
\end{ex}
Let $\X$ be a chainable relational structure and
$$
\L _\X= \Big\{ \la X, \vartriangleleft \ra : \;\vartriangleleft \mbox{ is a linear order on }X \mbox{ and chains }\X \Big\}.
$$
The following statement follows from Theorem 9 of  \cite{Gib},
which is a modification of the description of indiscernible sequences obtained independently by Frasnay in \cite{Fras}
and by Hodges, Lachlan and Shelah in \cite{Hodg1} (see also \cite{Fra}, p.\ 378).
\begin{fac}[Gibson, Pouzet and Woodrow]\label{T8090}
If $\X$ is a chainable relational structure and $\BL \in \L _\X$, then one of the following conditions holds:
\begin{itemize}
\item[{\sc (i)}]  $\displaystyle\L _\X =\LO _X$,
\item[{\sc (ii)}] $\displaystyle\L _\X =\bigcup _{\BL =\BI +\BF}\Big\{ \BF + \BI ,\, \BI ^* +\BF ^*\Big\}$,\\[-3mm]
\item[{\sc (iii)}]
$\displaystyle
\exists n\in \o \;\;\L _\X =\bigcup _{{\BL =\BK +\M +\BH \,,\; |K\cup H| \leq n }\atop  \vartriangleleft _K \in \LO_K ,\; \vartriangleleft _H \in \LO_H}
\Big\{ \BK _{\vartriangleleft _K}+ \M +\BH _{\vartriangleleft _H},\; \BH _{\vartriangleleft _H}+ \M ^* +\BK _{\vartriangleleft _K}\Big\}.
$
\end{itemize}
\end{fac}
\paragraph{Linear orders. Indecomposable ordinals}
A linear order $\BL=\la L, < \ra$ is called:
{\it dense} iff for each $x,y\in L$ satisfying $x<y$ there is $z\in L$ such that $x<z<y$;
{\it scattered} iff it does not contain a dense suborder iff $\Q \not\hookrightarrow \BL$;
{\it additively indecomposable} iff for each decomposition
of $\BL$ into an initial and a final part, $\BL=\BI + \BF$, we have $\BL\hookrightarrow \BI$ or $\BL\hookrightarrow \BF$; then we write $\BL \in AI$.
\begin{fac}\label{T4400}
For a limit ordinal $\a$ the following conditions are equivalent:

(a) $\b + \g < \a$, for each $\b , \g <\a$ ($\a $ is  indecomposable),

(b) $\b +\g= \a  \land \g>0 \Rightarrow \g =\a$ ($\a $ is right indecomposable),

(c) $\,\b +\a =\a$, for each $\b <\a$,

(d) $\a = \o ^\delta$, for some ordinal $\delta >0$,

(e) $\a\hookrightarrow A$ or $\a\hookrightarrow B$, whenever $\a =A\;\dot{\cup}\; B$ ($\a$ is an indivisible structure),

(f) $\I _\a =\{ I\subset \a : \a \not\hookrightarrow I\}$ is an ideal in $P(\a )$.
\end{fac}
\dok
For (a) $\Leftrightarrow $ (d) $\Leftrightarrow $ (c) see \cite{Sier}, p.\ 282, 323. For (b) $\Leftrightarrow $ (d) see \cite{Rosen}, p.\ 176.
The equivalence (a) $\Leftrightarrow $ (e) is 6.8.1 of \cite{Fra}. (e) $\Leftrightarrow$ (f) is evident.
\hfill $\Box$
\paragraph{Well-quasi-orders and Laver's theorem}
If $\P =\la P,\leq\ra$ is a quasi-order and $p,q\in P$, then we will write $p \mid q$ iff $p$ and $q$ are {\it incomparable} (i.e.\ $\neg p\leq q$ and $\neg q\leq p$);
and $p<q$ iff $p\leq q$ and $\neg q\leq p$. A set $A\subset P$ is an {\it antichain} iff $p \mid q$, for different $p,q\in A$.
A sequence $\la p_n :n\in \o\ra$ in $P$ is a {\it decreasing $\o$-sequence} iff $p_0 >p_1 >\dots$.
The quasi-order $\P$ is a {\it well-quasi-order (wqo)} iff
$\P$ has no infinite antichain and no decreasing $\o$-sequence.
An element $p$ of $P$ is a {\it minimal element of $\P$} iff there is no $q\in P$ such that $q<p$.
\begin{fac}\label{T001}
If $\P =\la P,\leq\ra$ is a wqo and $\emptyset \neq A\subset P$, then $\A =\la A , \leq \upharpoonright A\ra$ is a wqo and has a minimal element.
\end{fac}
\dok
The first claim is evident. Assuming that there are no minimal elements in $\A$,
for each $p\in A$ there would be $q\in A$ such that $q<p$; so, there would be a decreasing $\o$-sequence in $\A$ and, hence, in $\P$, which is false.
\kdok
It is evident that $\la \LO ,\hookrightarrow \ra$ is a quasi-order.
Let $\s$-$\Scatt$ denote the class of all linear orders which can be presented as unions of $\leq \o$ scattered linear orders.
\begin{fac}[Laver\cite{Lav}]\label{T002}
The quasi-order $\la \s\mbox{-}\!\Scatt , \hookrightarrow \ra$ is a wqo.
\end{fac}
\paragraph{Partial orders and forcing}
If $\P =\la P, \leq \ra $ is a preorder, the elements $p$ and $q$ of $P$ are {\it incompatible},
we write $p\perp q$, iff there is no $r\in P$ such that $r\leq p,q$.
An element $p$ of $P$ is an {\it atom } iff each two elements $q,r\leq p$ are compatible;
${\mathbb P} $ is called {\it atomless} iff it has no atoms.
A set $D\subset P$ is called: {\it dense} iff for each $p\in P$ there is $q\in D$ such that $q\leq p$;
{\it open} iff $q\leq p\in D$ implies $q\in D$.
The {\it distributivity number} $\h (\P )$ of $\P$
is the minimal size of a family $\CD$ of open dense subsets of $\P$ such that $\bigcap \CD$ is not dense;
in particular, $\h :=\h ((P(\o )/\Fin )^+)$ and $\o_1 \leq \h \leq \c$.
If $\kappa $ is a cardinal, ${\mathbb P} $ is called $\kappa ${\it -closed} iff for each
$\gamma <\kappa $ each sequence $\langle  p_\alpha :\alpha <\gamma\rangle $ in $P$, such that $\alpha <\beta \Rightarrow p_{\beta}\leq p_\alpha $,
has a lower bound in $P$. $\omega _1$-closed preorders are called {\it $\sigma$-closed}.

A preorder ${\mathbb P} =\langle  P , \leq \rangle $ is called
{\it separative} iff for each $p,q\in P$ satisfying $p\not\leq q$ there is $r\leq p$ such that $r \perp q$.
The {\it separative modification} of ${\mathbb P}$
is the separative preorder $\mathop{\rm sm}\nolimits ({\mathbb P} )=\langle  P , \leq ^*\rangle $, where
$p\leq ^* q \Leftrightarrow \forall r\leq p \; \exists s \leq r \; s\leq q $.
The {\it separative quotient} $\sq (\P)$ of ${\mathbb P}$
is the antisymmetric quotient of $\sm (\P)$ and $\ro (\sq (\P))$ is the {\it Boolean completion} of $\P$.
Preorders $\P$ and $\Q$ are called {\it forcing equivalent}, in notation $\P \equiv_{forc}\Q$, iff they produce the same generic extensions.
It is a standard fact that
${\mathbb P}\equiv_{forc}\mathop{\rm sm}\nolimits ({\mathbb P})\equiv_{forc}\mathop{\rm sq}\nolimits  ({\mathbb P})\equiv_{forc}\ro(\mathop{\rm sq}\nolimits  ({\mathbb P}))$
(see \cite{Kun}).

If $\Q $ is a preorder, a mapping
$f: \P \rightarrow \Q$ is a {\it complete embedding}, in notation $f:\P \hookrightarrow _c \Q$  iff
(ce1) $p_1 \leq  p_2 \Rightarrow f(p _1)\leq  f(p_2)$,
(ce2) $p_1 \perp p_2 \Leftrightarrow f(p _1)\perp f(p_2)$,
(ce3) $\forall q\in \Q \; \exists p\in \P \; \forall p' \leq  p \;\; f(p') \not\perp q$.
If $\P \hookrightarrow _c \Q$, then $\Q\equiv _{forc}\P \ast \pi $, where
$\pi $ is a $\P$-name for a preorder and $\P \ast \pi $ denotes the two-step iteration (see \cite{Kun}).

\section{Case I: Constant structures. Copy-dense structures}\label{S3}
Here we consider Case ({\sc i}) from Fact \ref{T8090} in a more general context.
Preliminarily we note that $\la \k , < \ra =\la \k ,\in\ra$ is the natural linear order on a cardinal $\k$,
$\Sym (\k )$ will denote the set of all bijections from $\k$ to $\k$
and $L_\emptyset$  will denote the ``empty language";
so, $L_\emptyset$-formulas are the formulas containing only the equality symbol.

Let $L$ be a relational language and $\Mod _L (\k )$ the set of $L$-structures with domain $\k$.
Clearly $\la \Mod _L (\k ), \hookrightarrow\ra$ is preorder,
its antisymmetric quotient is the partial order $\Mod _L (\k )/\!\rightleftarrows$,
and in between we have the quotient $\Mod _L (\k )/\!\cong$, which is a preorder.
We consider the following properties of a structure $\X\in\Mod _L (\k )$:

- $\forall \Y \in  \Mod _L (\k )\;\; ( \Y \hookrightarrow \X \Rightarrow  \Y \rightleftarrows \X )$  ($\X$ is {\it minimal in $\Mod _L (\k )$}),

- $\forall \Y \in  \Mod _L (\k )\;\; ( \Y \hookrightarrow \X \Rightarrow  \Y \,\cong\, \X )$  ($[\X]$ is {\it minimal in $\Mod _L (\k )/\!\cong$}),

- $\forall \Y \in  \Mod _L (\k )\;\; ( \Y \rightleftarrows \,\X \Rightarrow  \Y \cong \,\X )$   ($\X$ is a {\it  CSB structure} in $\Mod _L (\k )$).

\noindent
Of course, CSB is an abbreviation for Cantor-Schr\"{o}der-Bernstein.
\begin{prop}\label{T009}
For a structure $\X \in \Mod _L (\k )$, where $\k \geq \o$, we have\\[-8mm]
\begin{itemize}\itemsep=-2mm
\item[\rm (a)] $\X$ is a CSB structure iff  $\,\P (\X )=\P (\X )\!\uparrow \;:=\{ A\subset \k : \exists C\in \P (\X )\; C\subset A\}$;\\[-4mm]
\item[\rm (b)] The following conditions are equivalent {\rm (constant structures)}:\\[-8mm]
\begin{itemize}\itemsep=-1mm
\item[\rm (i)] $\X$ is chainable and $\L _\X =\LO _\k$ (i.e.\ in Fact \ref{T8090} we have Case {\rm ({\sc i})}),
\item[\rm (ii)] $\Aut  (\X )=\Sym (\k ) $ (such structures are called {\rm constant} by Fra\"{\i}ss\'{e}),
\item[\rm (iii)] Each partial bijection from $\k$ to $\k$ is a partial automorphism of $\X$,
\item[\rm (iv)] $\X$ is $\Sigma _0$-definable on its domain $\k$ by $L_\emptyset$-formulas;
\end{itemize}
\item[\rm (c)] The following conditions are equivalent {\rm (copy-maximal structures)}:\\[-8mm]
\begin{itemize}\itemsep=-1mm
\item[\rm (i)] $\X$ is chainable and $\la \k , < \ra\in \L _\X$,
\item[\rm (ii)] $\P (\X )= [\k]^{\k}$,
\item[\rm (iii)] $\P (\X )$ is a dense set in the poset $\la [\k]^{\k}, \subset \ra $ and $\,\X$ is a CSB structure,
\item[\rm (iv)] $\X$ is minimal in $\Mod _L (\k )/\!\cong$;
\end{itemize}
\item[\rm (d)] The following conditions are equivalent {\rm (copy-dense structures)}:\\[-8mm]
\begin{itemize}\itemsep=-1mm
\item[\rm (i)] $\P (\X )$ is a dense set in the poset $\la [\k]^{\k}, \subset \ra $,
\item[\rm (ii)] $\X$ is minimal in $\Mod _L (\k )$;
\end{itemize}
\item[\rm (e)] $\X$ is constant $\Rightarrow $ $\X$ is copy-maximal $\Rightarrow $ $\X$ is copy-dense.
\end{itemize}
\end{prop}
\dok
(a) If $\X$ is a CSB structure, $C\in \P (\X )$ and $C\subset A\subset \k$, then $\A\in \Mod _L (\k )$ and $\A \rightleftarrows \X$,
which implies that $\A\cong \X$ and, hence, $A\in \P (\X)$.
Conversely, if $\,\P (\X )=\P (\X )\!\uparrow$, $f:\X \hookrightarrow \Y $ and $g:\Y \hookrightarrow \X$,
then $\P (\X)\ni g[f[\k]]\subset g[Y]\subset \k$,
which implies that $g[Y] \in \P (\X)$ and, hence, $\Y \cong \X$.

(b) For the equivalence (i) $\Leftrightarrow$ (ii) see 9.5.2 of \cite{Fra},
while (ii) $\Leftrightarrow$ (iii) is 9.4.2(2) of \cite{Fra}.
If (iii) holds, then  the structure $\X$ is freely interpretable in the $L_\emptyset$-structure $\k$ (see 9.2.1 of \cite{Fra});
that is, (iv) is true (see \cite{Fra}, p.\ 6).
If  $\f (v_0,\dots ,v_{n-1})$ is an $L_\emptyset$-formula without quantifiers, $f\in \Sym (\k)$ and $\bar x \in \k ^n$,
then, clearly, $\k \models \f [\bar x]$ iff $\k \models \f [f\bar x]$; thus (iv) implies (ii).

(c) If (i) holds, then $[\k]^{\k}=\P (\la \k , < \ra)\subset \P (\X)$ (see Proposition \ref{T015}(a)) so we have (ii).
The implication (ii) $\Rightarrow$ (i) follows from the work of Gibson, Pouzet and Woodrow, see Theorem 25 of \cite{Gib}.
The equivalence (ii) $\Leftrightarrow$ (iii) follows from (a);  (ii) $\Leftrightarrow$ (iv) is evident.

(d) If (i) holds and $f: \Y \hookrightarrow \X$, then $f[Y]\in [\k]^\k$,
there is $A\in \P(\X)$ such that $A\subset f[Y]$ and, hence, $\X\cong \A \hookrightarrow \Y$.
Thus $\Y \rightleftarrows \X$ and (ii) is true.
Conversely, if (ii) holds and $A\in [\k ]^\k$, then $\A \hookrightarrow \X$
and by (ii) there is an embedding $f:\X \hookrightarrow \A$.
So, $\P (\X )\ni f[X]\subset A$ and (i) is true.

(e) is true because (b)(i) implies (c)(i) and (c)(ii) implies (d)(i).
\kdok
We recall that if $\k\geq 2$ and $\l \geq \o$ are cardinals,
then the {\it collapsing algebra} $\Col (\l ,\k)$ is the Boolean completion  of the reversed tree $\la {}^{<\l }\k ,\supset\ra$.
For $\k\geq \o$ the partial order $(P(\k )/[\k ]^{<\k })^+$ will be denoted by $\CP _\k $.
It is known that $\CP _\k$ is atomless, homogeneous, of size $2^\k$ and $\k ^{++}\leq \cc (\CP _\k)\leq (2^\k) ^+$ (see \cite{Balc1}, p.\ 372).
\begin{te}\label{T010}
If $\k \geq \o$ is a cardinal and $\X \in \Mod _L (\k )$ is a copy-dense (in particular, a constant or a copy-maximal) structure, then we have

(a) $\B _\X \cong \ro (\CP _\k)$;

(b) $\B _\X \cong \ro (P(\o )/\Fin)$, if $\k =\o$;

(c) If, in addition, $2^\k =\k ^+$ and $2^{\cf(\k )} =\cf (\k )^+$, then
\begin{equation}\label{EQ003}
\B _\X \cong
               \left\{\begin{array}{ll}
                 \Col (\o ,2^\k),    & \mbox{if }\cf (\k )>\o, \\
                 \Col (\o _1 ,2^\k), & \mbox{if }\cf (\k )=\o.
               \end{array}
               \right.
\end{equation}
\end{te}
\dok
(a) By the assumption, $\P (\X )$ is a dense set in the partial order $\la [\k]^{\k}, \subset \ra$,
which implies that $\P (\X )\equiv_{forc}\la [\k]^{\k}, \subset \ra$.
So, since $\sq (\la [\k]^{\k}, \subset \ra) = \CP _\k $,
we obtain  $\B _\X =\ro (\sq (\P (\X )))\cong \ro (\sq (\la [\k]^{\k}, \subset \ra))=\ro (\CP _\k)$.

(b) follows from (a) and the equality $\CP _\o= (P(\o )/\Fin)^+$.

(c) Under the assumptions, (\ref{EQ003}) holds for the algebra $\ro (\CP _\k )$ instead of $\B _\X$ (see \cite{Balc1}, p.\ 380).
Thus the statement follows from (a).\footnote{
More results similar to (\ref{EQ003}) can be obtained from the following statements.
\begin{fac}\label{T4130}
 Let $\l\geq \o$ be a regular cardinal and $\P$ a separative $\l$-closed preorder of size $\k=\k^{<\l}$.

(a) If $\k >\l$ and $1_\P\Vdash |\check{\k}|= \check{\l}$, then $\ro (\sq (\P ))\cong \Col (\l ,\k)$;

(b) If $\k =\l$ and $\P$ is atomless, then $\ro (\sq (\P ))\cong \Col (\l ,\l)$ (see \cite{KCol}).
\end{fac}
\begin{fac}\label{T209}
If $\k$ is an infinite cardinal, then for the partial order $\CP _\k$ we have

(a)  Forcing by $\CP _\o$ collapses $\c$ to $\h$ (Balcar, Pelant and Simon, \cite{BPS});

(b)  If $\k >\cf (\k ) =\o$, then $\Col(\o _1, \k ^\o)\hookrightarrow _c \ro (\CP _\k)$ (Kojman, Shelah \cite{Koj});

(c)  If $\k > 2^{\cf (\k )}> \cf (\k ) >\o$, then  $\Col(\o , \k ^+)\hookrightarrow _c \ro (\CP _\k)$ (Shelah, \cite{She});

(d)  If $\k = \cf (\k )>\o$, then $\CP _\k $ collapses each $\lambda <\cc (\CP _\k)$ to $\o$ (Shelah, \cite{She}).
\end{fac}
}
\hfill $\Box$
\begin{ex}\label{EX000}
Copy-dense structures which are neither chainable nor copy-maximal. \rm
The disjoint union of linear orders $\X =\dot{\bigcup} _{n\in \o}\BL_n$, where $\BL_n \cong\o _n$, for $n\in \o$, is not 2-monomorphic,
because it contains both comparable and incomparable pairs; by Proposition \ref{T009}(c) $\X$ is not copy-maximal
and, clearly, $|X|=\aleph _\o$. If $A\in [X]^{\aleph _\o}$, then $\sum _{n\in \o}|A\cap L_n|=\aleph _\o$
and an easy recursion gives an increasing sequence $\la n_k :k\in \o\ra$ in $\o$ such that $|A\cap L_{n_k}|\geq \aleph _k$, for each $k\in \o$.
Consequently, the set $A$ contains a copy of $\X$ and $\X$ is a copy-dense structure.
By Theorem \ref{T010} we have $\B _\X \cong \ro (\CP _{\aleph _\o})$;
if, in addition, $\aleph _\o ^\o =2^{\aleph _\o}$, then by Facts \ref{T4130} and \ref{T209}(b) we have $\B _\X \cong \Col (\o _1, 2^{\aleph _\o}) $.
\end{ex}
\begin{rem}\label{R001}\rm
If $\X \in \Mod _{L}(\k)$ is a chainable structure and there is a linear order $\BL \in \L _\X$ which is copy-dense (that is minimal in $\Mod _{L_b}(\k)$),
then $\X$ is copy-dense as well, because by Proposition \ref{T015}(a) we have $\P (\BL) \subset \P (\X)$.
We note that the set of $\c$-sized suborders of $\BR$ does not have minimal elements (Dushnik and Miller, \cite{Dus}; see also \cite{Fra}, p.\ 150)
and under PFA all $\aleph _1$-dense suborders of $\BR$ are isomorphic and minimal (Baumgartner, \cite{Baum}).
So if a structure $\X$ is chained by one of them, then by Theorem \ref{T010} we have $\B_\X \cong \ro (\CP _{\o _1})\cong \Col (\o ,\c)$, because under PFA we have $2^{\o _1}=\o _2 =\c$.
The same holds if instead of an $\aleph _1$-dense set we take a Countryman line
(see Moore's article \cite{Moo}).
\end{rem}
\begin{ex}\label{EX001}
Chainable copy-dense structures which are not copy-maximal. \rm
The sum of linear orders $\BL =\sum _{n\in \o}\BL_n$, where $\BL_n \cong\o _n ^*$, for $n\in \o$, is not CSB (e.g.\ $1+\BL \rightleftarrows \BL$, but $1+\BL \not\cong \BL$)
and, by Proposition \ref{T009}(c), $\BL$ is not copy-maximal.
Clearly, $|L|=\aleph _\o$ and as in Example \ref{EX000} we easily show that $\BL$ is copy-dense.
If $\X =\la L, D_{\f _c}\ra$ is the cyclic relation defined on $\BL$,
then by Proposition \ref{T015}(a) we have $\P (\BL) \subset \P (\X)$ and, hence, $\X$ is copy-dense as well.
By Theorem \ref{T010} we have $\B _\X \cong \ro (\CP _{\aleph _\o})$.
\end{ex}
\begin{ex}\label{EX002}
CSB structure which is not chainable. \rm
The disjoint union of linear orders $\X =\BL _0 \,\dot{\cup}\, \BL _1$, where $\BL_0 \cong\o $ and $\BL_1 \cong\o ^*$ is not 2-monomorphic.
It is evident that $\P (\X)\!\uparrow =\P (\X)$ and, by Proposition \ref{T009}(a) $\X$ is a CSB structure.
\end{ex}
\begin{rem}\label{R002}\rm
Concerning CSB structures we note that countable CSB linear orders are scattered,
but  there is a dense suborder $\BE$ of the real line which is embedding-rigid; that is, $\Emb (\BE)=\{ \id _E\}$
(Dushnik and Miller \cite{Dus}; see also \cite{Rosen}, p.\ 147);
thus, $\BE$ is a CSB structure, because $\P (\BE)\!\uparrow =\P (\BE)=\{ E\}$.
For scattered linear orders we have:
A scattered linear order is CSB iff it is a finite sum of
well orders, their inverses, l.o.'s of the form $\o ^\theta \o ^* +\o ^{\d }$, where $1\leq\theta <\d$ are ordinals, and their inverses
(Laflamme, Pouzet and Woodrow \cite{Laf}).
\end{rem}
\section{Question \ref{Q1}: dense embedding of copies of chains}\label{S4}
In this section, using the following basic statement,
we show that Question \ref{Q1} has the answer `Yes' for monomorphic structures belonging to a large class.
\begin{prop}\label{T015}
If $\X$ is a chainable structure and $\BL \in \L _\X$, then

(a) $\Emb (\BL )\subset \Emb (\X )$ and $\P (\BL )\subset \P (\X )$;

(b) If $\P (\BL)$ is a dense subset of $\P (\X)$, then $\B_\X \cong \B_\BL$.
\end{prop}
\dok
(a) If  $f\in \Emb (\BL )$, then $f$ is a partial automorphism of $\BL$.
Since $\BL \in \L _\X$  we have $\Pa (\BL)\subset \Pa (\X )$;
thus $f\in \Pa (\X )$ and, since $\dom (f)=X$,  $f\in\Emb (\X )$.
So, we have $\P (\BL )=\{ f[X]: f\in \Emb (\BL )\}\subset \{ f[X]: f\in \Emb (\X )\}=\P (\X )$.

(b) If $\P (\BL)$ is dense in $\P (\X)$, then $\P (\BL)\equiv_{forc}\P (\X)$ and, by (\ref{EQ005}), $\B_\X \cong \B_\BL$.
\kdok
Let $\X $ be a chainable structure and $\BL =\la X ,<\ra\in \L _\X$.
If $f\in \Emb (\X)$,
then $f\times f:X^2 \rightarrow X^2$, where $(f\times f)(\la x,x'\ra)=\la f(x),f(x')\ra$;
instead of $(f\times f)^{-1}[<]$ we will write only $f^{-1}[<]$;
namely, $\la x, x'\ra \in f^{-1}[<]$ iff $f(x)<f(x')$.
\begin{lem}\label{T016}
If $\X $ is a chainable structure, $\BL =\la X,< \ra \in \L _\X$ and $f\in \Emb (\X)$, then
\begin{equation}\label{EQ012}
\BL _f :=\la X , f^{-1}[<]\ra \in \L _\X \;\;\mbox{ and }\;\;f:\BL _f \hookrightarrow \BL .
\end{equation}
\end{lem}
\dok
Since $f\in \Emb (\X)$ the mapping $f:X\rightarrow X$ is an injection.
For $x,x'\in X$ we have $\la x,x'\ra\in f^{-1}[<]$ iff $f(x)<f(x')$; thus, $f:\BL _f \hookrightarrow \BL$.

For a proof that $\BL _f \in \L _\X$
we show that the structure $\X=\la X,\bar \r \ra$ is $\Sigma _0$-definable in the linear order $\BL _f$.
Since $\BL \in \L _\X$
by Fact \ref{T8067}  for $i\in I$ there is an $L_b$-formula without quantifiers
$\f _i (v_1, \dots ,v_{n_i})=\f _i (\bar v)$ such that
\begin{equation}\label{EQ000}
\forall \bar{x}\in X^{n_i}\;\; ( \bar{x}\in \r _i \Leftrightarrow \BL \models \f _i [\bar{x}]) .
\end{equation}
So, for each $\bar{x}\in X^{n_i}$ we have
$\bar{x}\in \r _i $
iff $f\bar{x}\in \r _i $ (since $f\in \Emb (\X)$)
iff $\BL \models \f _i [f\bar{x}]$ (by (\ref{EQ000}))
iff $\BL _f \models \f _i [\bar{x}]$ (because $f:\BL _f \hookrightarrow \BL$).
Thus for each $i\in I$  (\ref{EQ000}) is true when we replace $\BL$ by $\BL _f$ and we are done.
\hfill $\Box$
\begin{lem}\label{T017}
If $\BL$ is a linear order and $f:\BL ^* \hookrightarrow\BL$, then $f:\BL  \hookrightarrow\BL^*$.
\end{lem}
\dok
Clearly, if $\BL _1$ and $\BL_2$ are linear orders and $f:\BL _1 \hookrightarrow \BL _2$,
then $f:\BL _1^* \hookrightarrow \BL _2^*$.
So, $f:\BL ^* \hookrightarrow\BL$ implies that $f:(\BL ^*)^*  \hookrightarrow\BL ^*$,
that is, $f:\BL  \hookrightarrow\BL^*$.
\kdok
\begin{te}\label{T018}
If $\X$ is a non-constant chainable structure,
then each of the following conditions implies that $\B _\X \cong \B _\BL$, for some $\BL\in \L _\X$

(a) There is a minimal element $\BL $ of the quasi-order $\la \L _\X ,\hookrightarrow\ra$; then $\B_\X\cong \B_\BL $;

(b) There is $\BL \in \L _\X$ satisfying {\sc (iii)}; then $\B_\X\cong \B_\BL $;

(c) There is $\BL \in  \L _\X \cap AI$; then $\B_\X\cong \B_\BL $;

(d) There is $\BL \in\L _\X \cap  \s\mbox{-}\!\Scatt $; then $\B_\X\cong \B_{\BL '}$, for some $\BL' \in\L _\X$.
\end{te}
\dok
The initial part of the proof in cases (a)--(c) is the same:
``If $\BL $ is as assumed, then by Proposition \ref{T015} we will have $\B _\X \cong \B _\BL$, if $\P (\BL )$ is dense in $\P (\X )$; that is,
\begin{equation}\label{EQ8138}
\forall C\in \P (\X) \;\; \exists A\in \P (\BL ) \;\; A\subset C.
\end{equation}
Let $C\in \P (\X )$ and $f\in \Emb (\X )$, where $f[X]=C$.
By Lemma \ref{T016} we have (\ref{EQ012})."
Here we show that the existence of $A\in \P (\BL )$ satisfying $A\subset C$  follows from
\begin{equation}\label{EQ013}
\BL \hookrightarrow \BL _f \;\;\mbox{ or }\;\;\BL ^* \hookrightarrow \BL _f.
\end{equation}
First, if $g: \BL \hookrightarrow \BL _f$,
then by (\ref{EQ012}) we have $f\circ g : \BL \hookrightarrow \BL$;
so $A=f[g[X]]\in \P (\BL )$ and $A\subset f[X]= C$.
Second, if $g: \BL  ^*\hookrightarrow  \BL _f$,
then $h:=f\circ g : \BL  ^*\hookrightarrow \BL$,
and, by Lemma \ref{T017}, $h :\BL \hookrightarrow\BL ^*$,
which implies that $h \circ h: \BL  \hookrightarrow \BL$.
So, for $A:= h[h[X]]$ we have $A\in \P (\BL )$
and $A\subset h[X]=f[g[X]]\subset f[X]= C$.

(a) Let $\BL$ be a minimal element of $\la \L _\X ,\hookrightarrow\ra$.
By (\ref{EQ012}) we have $\BL _f \hookrightarrow \BL$;
so, by the minimality of $\BL$ we have $\BL \hookrightarrow \BL _f$ and (\ref{EQ013}) is true.

(b) Let ({\sc iii}) hold for $\BL \in \L _\X$. By (\ref{EQ012}), $\BL _f\in \L _\X$ and we have two cases.

1. $\BL _f= \la K, \vartriangleleft _K \ra +\M + \la H, \vartriangleleft _H \ra$, where $\BL  =\BK +\M +\BH$ and $|K|,|H|<\o$.
Then $\BL \cong \BL _f $, which implies (\ref{EQ013}).

2. $\BL _f= \la H, \vartriangleleft _H \ra +\M ^* + \la K, \vartriangleleft _K \ra$, where $\BL  =\BK +\M +\BH$ and $|K|,|H|<\o$.
Then $\BL ^*  =\BH ^* +\M ^* +\BK ^*$ and, hence,  $\BL ^* \cong \BL _f$, which implies (\ref{EQ013}).

(c) Let $\BL \in  \L _\X \cap AI$. If ({\sc iii}) holds for $\BL$, we apply (b).
Otherwise ({\sc ii}) holds for $\BL$
and, since by  (\ref{EQ012}) we have $\BL _f \in \L _\X$, we have the following cases.

1. $\BL _f=\BF + \BI$, where $\BL =\BI +\BF$.
Then, since $\BL \in AI$, we have $\BL \hookrightarrow \BI$ or $\BL \hookrightarrow \BF$;
so, $\BL\hookrightarrow \BL _f$ and (\ref{EQ013}) is true.

2. $\BL _f=\BI ^* + \BF^*$, where $\BL =\BI +\BF$.
Then again $\BL \hookrightarrow \BI$ or $\BL \hookrightarrow \BF$,
and, hence,  $\BL ^* \hookrightarrow \BI^*$ or $\BL ^* \hookrightarrow \BF^*$,
which implies that $\BL  ^* \hookrightarrow \BL _f $ and (\ref{EQ013}) is true again.

(d) If $\BL \in \L _\X \cap  \s\mbox{-}\!\Scatt$,
then $L=\bigcup _{n\in \o}S_n$, where $S_n$, $n\in \o$, are scattered suborders of $\BL$.
If ({\sc iii}) holds for $\BL$, we apply (b).
Otherwise ({\sc ii}) holds for $\BL$
and $\L _\X \subset  \s\mbox{-}\!\Scatt$.
E.g.\  if  $\BL =\BI +\BF$,
then $\BF +\BI =\bigcup _{n\in \o}(S_n \cap F) + \bigcup _{n\in \o}(S_n \cap I)\in \s\mbox{-}\!\Scatt$,
because a suborder of a scattered order is scattered
and the sum of two scattered orders is scattered.
By Facts \ref{T001} and \ref{T002} there is a minimal element $\BL '$ in the quasi-order $\la \L _\X, \hookrightarrow\ra$
and the statement follows from (a).
\hfill $\Box$
\section{Structures chained by ordinals}\label{S5}
\begin{lem}\label{T8115}
If $\d \geq 1$ is an ordinal, $s\in \N$ and $\o ^\d s = \b +\g$, then $\g +\b \geq \o ^\d s$.
\end{lem}
\dok
If $\b <\o ^\d$, then $\o ^\d s = \b +\g _1 +\g_2$, where $\o ^\d  = \b +\g _1$ and, hence, $\g _1=\o ^\d$, and $\g _2= \o ^\d (s-1)$.
So $\g =\g_1 +\g _2=\o ^\d s$ and $\g +\b \geq \o ^\d s$.

If $\o ^\d \leq \b < \o ^\d s$, then dividing $\b$ by $\o ^\d$ we have
$\b =\o ^\d s_1 +\theta$, where $1\leq s_1 <s$ and $0\leq\theta< \o ^\d$.
If $\theta =0$, then by the uniqueness of the difference of ordinals
we have $\g = \o ^\d (s-s_1)$ and $\g +\b =\o ^\d s$.

If $\theta >0$, then $\o ^\d s_1 +\o ^\d (s-s_1)=\o ^\d s_1 + \theta +\g$
and, by the left cancelation law, $\o ^\d (s-s_1)=\theta +\g$.
Assuming that $\g <\o ^\d$ by Fact \ref{T4400}(a) we would have $\theta +\g <\o ^\d \leq \o ^\d (s-s_1)$
which is false; thus $\g \geq \o ^\d$, which gives $\g =\o ^\d +\g '$.
So $\o ^\d (s-s_1)=\theta +\o ^\d +\g '$
and, since by Fact \ref{T4400}(c) we have $\theta +\o ^\d=\o ^\d$,
we obtain $\o ^\d (s-s_1)=\o ^\d +\g '=\g$.
Thus $\g +\b = \o ^\d s + \theta > \o ^\d s$.
\hfill $\Box$
\begin{te}\label{T012}
If $\a =\o ^{\d _n }s_n + \dots + \o ^{ \d _0 }s_0 +k $ is an infinite ordinal
presented in the Cantor normal form, where $n, k\in \o$, $s_0,\dots ,s_n \in \N$,
and $0<\d _0<  \dots < \d_n$ are ordinals,  then

(a) $\sq \P (\a )\cong \prod _{i=0}^n \big(\sq \P (\o ^{\d _i})\big)^{s_i}$;

(b) $\min \{\g +\b : \a =\b +\g \}=\o ^{\d _n }s_n $;

(c) If $\X$ is a chainable structure and $\a \cong \BL \in \L _\X$, then
\begin{equation}\label{EQ008}
\B _\X \cong
               \left\{\begin{array}{cl}
                 \ro \Big(\CP _{|\a|}\Big),                                         &\mbox{if ({\sc i}) holds},\\
                 \ro \Big(\sq (\P (\o ^{\d _n}))^{s_n}\Big),                        & \mbox{if ({\sc ii}) holds} , \\
                 \ro \Big(\prod _{i=0}^n \big(\sq \P (\o ^{\d _i})\big)^{s_i}\Big), & \mbox{if ({\sc iii}) holds}.
               \end{array}
               \right.
\end{equation}
\end{te}
\dok
(a) is Theorem 3.2 of \cite{Kord1} (and has a direct proof).

(b) ($\leq$) For $\b =\o ^{\d _n}s_n $ and $\g=\o ^{\d _{n-1} }s_{n-1} + \dots + \o ^{ \d _0 }s_0 +k$
we have $\a = \b +\g $
and, by the Cantor normal form theorem, $\g  < \o ^{\d _n }$,
which by Fact \ref{T4400}(c) implies that $\g +\o ^{\d _n }=\o ^{\d _n }$
and, hence, $\g +\b =\g +\o ^{\d _n } + \o ^{\d _n  }(s_n -1)=\o ^{\d _n } + \o ^{\d _n  }(s_n -1)=\o ^{\d _n }s_n$.

($\geq$) Let $\a = \b +\g$. If $\b \geq \o ^{\d _n }s_n$, then $\g +\b \geq \o ^{\d _n }s_n$ and we are done.
If $\b <\o ^{\d _n }s_n$,
then, since $\a = \b +\g$, we have $\g =\g _1 +\g _2$, where $\b +\g _1=\o ^{\d _n }s_n$;
by Lemma \ref{T8115} we have  $\g _1 +\b\geq \o ^{\d _n }s_n$
and, since $\g \geq \g _1$, we have $\g  +\b\geq \o ^{\d _n }s_n$.

(c) If ({\sc i}) holds the statement follows from Proposition \ref{T009} and Theorem \ref{T010}.

Let ({\sc ii}) hold. Taking a bijection from $X$ to $\a$
w.l.o.g.\ we can suppose that $X=\a$
and that $\BL =\la \a, <\ra$, where $<$ is the natural order on $\a$.
Since $\BL$ satisfies ({\sc ii}) of Fact \ref{T8090} we have
$\L _\X =\bigcup _{\a =\BI +\BF}\{ \BF + \BI ,\, \BI ^* +\BF ^*\}$.
Here $\a =\BI +\BF$ means that there are ordinals $\b$ and $\g$ such that $\BI\cong \b$, $\BF \cong \g$ and $\a=\b +\g$,
which means that $\a =\b +(\a-\b)=[0,\b) +[\b, \a)$ and $\g \cong [\b, \a)$.
Thus $\BF +\BI$ denotes the linear order on the set $\a$ of the form $[\b, \a)+[0,\b)$ isomorphic to the ordinal $\g +\b$,
while $\BI^* +\BF^*\cong \b ^* +\g ^*$.
So, if $\ord$ denotes the function which to each well order $\BL$ adjoins the ordinal $\ord(\BL)\cong \BL$, we have
\begin{equation}\label{EQ007}\textstyle
\ord [\L _\X] =\bigcup _{\a =\b +\g}\{ \g + \b ,\, \b ^* +\g ^*\}.
\end{equation}
By (a) there is a linear order $\vartriangleleft $ on the set $\a$
such that $\BL _m :=\la \a, \vartriangleleft \ra\cong \o ^{\d_n }s_n $ and $\BL _m \in \L_\X$.
In order to apply Theorem \ref{T018}(a) we show that $\BL _m$ ia a minimal element of the quasi order $\la \L _\X , \hookrightarrow\ra$.
So, if $\BL '\in \L _\X$ and $\BL '\hookrightarrow \BL _m$,
then, since $\BL _m$ is a well order and $\a \geq \o$,
for each partition $\a =\b +\g$ we have $\otp (\BL ')\not\cong \b ^* +\g ^*$.
Thus, by (\ref{EQ007}) we have $\otp (\BL ')\cong \g +\b$, where $\a =\b +\g$,
which by (a) implies that $\o ^{\d_n }s_n \leq \g +\b$
and, hence, $\BL _m \hookrightarrow \BL'$.
So, by Theorem \ref{T018}(a) we have $\B _\X \cong \B _{\BL _m}=\ro(\sq(\P (\o ^{\d_n }s_n)))$.
By (a) we have $\sq(\P (\o ^{\d_n }s_n))\cong \sq (\P (\o ^{\d _n}))^{s_n}$.

If ({\sc iii}) holds, the statement follows from (a) and Theorem \ref{T018}(b).
\hfill $\Box$
\section{Countable monomorphic structures}\label{S6}
Here by Theorem \ref{T018}(d) we have $\B _\X \cong \B _\BL$, for some $\BL \in \L _\X$.
We recall that sh$(\S )$ denotes the cardinality of the ground model ${\mathfrak c}$ in the Sacks extension.
\begin{te}\label{T000}
If $\X$ is a countable chainable structure, then

(a) There is a l.\ o.\ $\BL \in \L _\X $ such that $\B_\X \cong \B _\BL$;

(b) If $\X$ is not a constant structure, then $\L _\X \subset \Scatt$ or $\L _\X \cap \Scatt=\emptyset$.
\end{te}
\dok
(a) Countable linear orders are $\s$-scattered and we apply Theorem \ref{T018}(d).

(b) If $\X$ is non-constant and chained by a linear order $\BL$,
then regarding Fact \ref{T8090}, for the set $\L_\X$ we have ({\sc ii}) or ({\sc iii}).
It is evident that $\BL \in \Scatt$ implies that $\L _\X \subset \Scatt$ in both cases.
\hfill $\Box$
\begin{te}\label{T011}
Let $\X$ be a non-constant countable chainable structure.
\begin{itemize}
\item[\rm (a)] If $\L _\X \subset \Scatt$, the algebra $\B _\X $ contains a dense $\s$-closed atomless subset.\\
               If CH holds, then $\B_\X \cong \ro (P(\o )/\Fin)$.
\item[\rm (b)] Otherwise, $\B _\X\cong \ro (\S \ast \pi)$, where $\pi$ is an $\S$-name for a $\s$-closed,
               separative and atomless poset of size $\c$.\\
               If $\,\sh (\S )=\aleph _1$ or PFA holds,
               then  $\B _\X\cong \ro (\S \ast \dot{(P(\o )/\Fin)^+})$,
               where $\dot{(P(\o )/\Fin)^+}$ is an $\S$-name for $(P(\o )/\Fin)^+$ in the Sacks extension.
\end{itemize}
\end{te}
\dok
 By Theorem \ref{T000}(a) there is a linear order $\BL \in \L _\X $ such that $\B_\X \cong \B _{\BL }$.

(a) If $\L _\X \subset \Scatt$, then $\BL  \in \Scatt$ and $\B _\X \cong \ro(\sq (\P (\BL )))$.
By the main result of \cite{Kscatt} the poset $\sq (\P(\BL  ))$ is $\s$-closed atomless of size $\c$.
Under CH, by Fact \ref{T4130} we have $\B _{\BL }\cong \Coll (\o _1,\o _1)\cong \ro (P(\o )/\Fin)$
and, hence, $\B_\X \cong \ro (P(\o )/\Fin)$.

(b) By Theorem \ref{T000}(b) $\BL $ is a non-scattered linear order. By the main results of \cite{KurTod}
the posets $\P (\BL )$ and, hence, $\sq(\P (\BL ))$ and $\B _{\BL }$, are forcing equivalent to $\S \ast \pi$, where $\pi$ is as above,
and all the statements of claim (b) of the present theorem are true for the algebra $\B _{\BL '}$;
thus, they are true for the algebra $\B _\X$ as well.
\kdok
We note that, regarding Theorem \ref{T011}(b), the equality $\sh (\S )=\aleph _1$ follows from CH and, more generally, from ${\mathfrak b}=\aleph _1$,
by a result of Simon \cite{Sim}. Then  we have $\B _\X\cong \ro (\S \ast \dot{\Col(\o _1,\o _1)^+})$,
while  $\B _\X\cong \ro (\S \ast \dot{\Col(\o _2,\o _2)^+})$, under the PFA.

So, under CH, for each non-constant countable chainable structure $\X$ we have
$$
\B _\X \cong
               \left\{\begin{array}{cl}
                 \ro (P(\o )/\Fin),                   & \mbox{if $\X$ is chained by a scattered l.\ o.,} , \\
                 \ro (\S \ast \dot{(P(\o )/\Fin)^+}), & \mbox{if $\X$ is chained by a non-scattered l.\ o.}
               \end{array}
               \right.
$$
For countable structures chained by well orders we can say more than the ZFC part of Theorem \ref{T011}(a) states.
We recall notation from \cite{Kord}.
If $\B$ is Boolean algebra, $\rp (\B )$ denotes the reduced power $\B ^\o / \Fin$,
and its iterations $\rp ^n(\B )$, for $n\in \o$, are defined by $\rp ^0(\B) =\B $ and $\rp ^{n+1}(\B )= \rp (\rp ^n (\B ))$.
For an ordinal $\d >0$, the set $\I_{\o ^\d}:=\{ I\subset \o ^\d :\o ^\d \not\hookrightarrow I \}\subset P(\o ^\d)$ is an ideal
and  $P(\o ^\d)/ \I_{\o ^\d}$ is the corresponding quotient algebra.
\begin{te}\rm\label{T8117}
If $\X$ is a countable relational structure chained by a countable ordinal
$\a =\o ^{\g _n +r_n }s_n + \dots + \o ^{ \g _0 +r_0 }s_0 +k $ presented in the Cantor normal form,
where $k\in \o$, $r_i \in \o$, $s_i \in \N$, $\g _i \in \Lim \cup \{ 1 \}$
and $\g _n +r_n > \dots > \g _0 +r_0$, then

(a) Regarding the cases from Fact \ref{T8090} we have
\begin{equation}\label{EQ014}
\B _\X \cong
               \left\{\begin{array}{cl}
                 \ro \Big(P(\o )/\Fin\Big),                                                                   &\mbox{if ({\sc i}) holds},\\
                 \ro \Big( ( ( \rp ^{r_n}( P(\o ^{\g _n} )/ \I _{\o ^{\g _n} }))^+ )^{s_n}\Big),              & \mbox{if ({\sc ii}) holds} , \\
                 \ro \Big(\prod _{i=0}^n (( \rp ^{r_i}( P(\o ^{\g _i} )/ \I _{\o ^{\g _i} }))^+ )^{s_i}\Big), & \mbox{if ({\sc iii}) holds};
               \end{array}
               \right.
\end{equation}

{\rm (b)} $\B _\X\cong \ro ((P(\o )/\Fin )^+ \ast \pi)$, where $[\o ] \Vdash$ ``$\pi$ is $\s$-closed and separative";

{\rm (c)} $\B _\X\cong \ro ((P(\o )/\Fin )^+)$, if, in addition, ${\mathfrak h}=\o _1$.
\end{te}
\dok
 The first claim from (a) follows from Theorem \ref{T010}.
For $\a$ in the given form by Theorem 3.1 of \cite{Kord} we have $\sq (\P (\a ))\cong \prod _{i=0}^n (( \rp ^{r_i}( P(\o ^{\g _i} )/ \I _{\o ^{\g _i} }))^+ )^{s_i}$.
So, the second claim follows from Theorem  \ref{T012}(c)
and the third from Theorem \ref{T018}(b).

By Theorems 5.3 and 5.4 of \cite{Kord} $\ro (\sq (\P (\a )))\cong \ro ((P(\o )/\Fin )^+ \ast \pi)$, where $\pi$ is a $(P(\o )/\Fin )^+$-name for a $\s$-closed separative poset
and if, in addition, ${\mathfrak h}=\o _1$, then $\ro (\sq (\P (\a )))\cong \ro ((P(\o )/\Fin )^+)$.
So statements (b) and (c) follow from Theorem \ref{T000}(a).
\kdok
Analysis of the algebras from (\ref{EQ014}) in several models of ZFC is given in \cite{Kord}.
\section{Uncountable monomorphic structures}\label{S7}
\paragraph{Structures chained by uncountable ordinals}
Here we consider a relational structure $\X$ chained by an uncountable ordinal
\begin{equation}\label{EQ009}
\a =\o ^{\d _n }s_n + \dots + \o ^{ \d _0 }s_0 +k ,
\end{equation}
where $n, k\in \o$, $s_0,\dots ,s_n \in \N$,
and $0<\d _0<  \dots < \d_n$ are ordinals.
By Theorem \ref{T012}(c) the algebra $\B _\X$ is isomorphic to the Boolean completion
of a finite direct product of posets of the form $\sq \P (\o ^{\d _i})$
and their ZFC properties are described in the following statement (Theorem 4.10 from \cite{Kord1}).
\begin{te}\label{T207}
If $\d>0$ is an ordinal and $\cf (\d )=\k $, then for $\P (\o ^\d)$ we have\\[-6mm]
\begin{itemize}\itemsep=-1.5mm
\item[\rm (A)] If $\d$ is a successor ordinal or $\cf (\d )=\o$ (that is, $\k \leq \o$), \\
               then $\CP _\o \hookrightarrow _c \sq \P (\o ^\d)$ and $\sq \P (\o ^\d)$ is $\sigma$-closed;
\item[\rm (B)] If $\d = \theta + \k$, where $\ord \ni\theta\geq \k >\cf (\theta )=\o$\\
               and $\theta =\lim _{n\rightarrow \o}\d _n$, where $\cf (\d _n)=\k$, for all $n\in \o$,\\
               then $\CP _\o \hookrightarrow _c \sq \P (\o ^\d)$ and $\sq \P (\o ^\d )$ is $\sigma$-closed;
\item[\rm (C)] If $\d = \theta + \k$, where $\ord \ni \theta \geq \k >\cf (\theta )=:\lambda > \o$\\
               and $\theta =\lim _{\xi\rightarrow \,\lambda }\d _\xi $, where $\,\cf (\d _\xi)=\k$, for all $\xi\in \lambda $,\\
               then $\CP _\lambda \hookrightarrow _c \sq \P (\o ^\d )$;
\item[\rm (D)] If $\d =\theta +\k$, where $\ord \ni \theta\geq \cf (\theta )\geq \k> \o$ or $\theta =0$,\\
               then $\CP _\k \hookrightarrow _c \sq \P (\o ^\d )$;
\item[\rm (E)] If $\d =\lim _{\xi \rightarrow \k} \d _\xi $, where $\cf (\d _\xi)=\k >\o$, for all $\xi< \k$,\\
               then $\CP _\k \hookrightarrow _c \sq \P (\o ^\d)$.
\end{itemize}
\end{te}
For example, if $|\d |=\o _2$, then
(A) holds for $\o _2 +1$ and $\o _2 +\o$,
(B) for $\o _2 \o +\o _2$,
(C) for $\o _2 \o _1 +\o _2$,
(D) for $\o _2 $ and $\o _2 +\o _2$,
and (E) holds for $\o _2 \o _2 $.
The following statement (Theorem 4.11 from \cite{Kord1})
describes the ZFC properties of the algebra $\B_\a=\ro(\sq (\P (\a )))$
and, together with Theorem \ref{T012}(c), of the algebra $\B_\X$ as well.
\begin{te}\label{T226}
If $\a$ is an uncountable ordinal given by (\ref{EQ009}), then

{\rm (i)} If $\d _i$ satisfies (A) or (B), for each $i\leq n$, then the partial order $\sq (\P (\a ))$ is $\s$-closed
and  ${\CP _\o }^k\hookrightarrow_c \sq (\P (\a ))$, where $k={\sum _{i=0}^n s_i}$;

{\rm (ii)} Otherwise, $\CP _\l\hookrightarrow _c\sq \P (\a )$, for some regular cardinal $\l >\o$ and collapses $\o _2$ to $\o$.
If, in addition, $\cc (\CP _\l)=(2^{|\a|})^+$, then $\B_\a \cong \Col (\o ,2^{|\a|})$.
\end{te}
\begin{ex}\label{EX003}
$\BL$ chains $\X$, but $\P (\BL)$ is not a dense subset of $\,\P (\X)$. \rm
Let $\BL=\la X,<\ra$ be a linear order and let $\X=\la X, D_{\f _c}\ra$,
where $D_{\f _c}$ is the cyclic relation defined on $\BL$ (see Example \ref{EX8008}).
It is easy to check that $\L _{\X } =\bigcup _{\BL =\BI +\BF}\{ \BF + \BI ,\, \BI ^* +\BF ^*\}$; that is, we have Case ({\sc ii}) in Fact \ref{T8090}.

Let $\BL $ be the ordinal $\a=\o _2\o +\o _2 =\o ^{\o _2 +1}+\o ^{\o_2}=\o ^{\d _1}+ \o ^{\d _0}$.
Then the exponent $\d _1$ satisfies (A) of Theorem \ref{T207},
the partial order $\sq (\P (\o ^{\d _1}))$ is $\s$-closed
and, by Theorem \ref{T012}(c), $\B _\X \cong \ro(\sq (\P (\o ^{\d _1})))$.
Thus the algebra $\B _\X$ is $\o$-distributive.
On the other hand, the exponent $\d _0$ satisfies (D) of Theorem \ref{T207},
and by Theorem \ref{T226} the partial order $\P (\BL)$ collapses $\o _2$ to $\o$.
Thus, $\P (\BL)$ is not a dense subset of the poset $\P (\X)$ and $\B _\X \not\cong \B_\BL$.
\end{ex}
Assuming more than ZFC (i.e.\ under additional assumptions concerning cardinal arithmetic)
we obtain several consistent characterizations of the algebras $\B _\X$.
\begin{ex}\label{EX004}\rm
By Theorem 4.12 of \cite{Kord1}, if the equalities $\h =\o _1 $ and $\c =\o _2= 2^{\o _1}$ hold, then for each ordinal $\a \in [\o , \c)$
we have $\B _\a \cong \Col (\o _1 ,\c)$  or $\B _\a \cong \Col (\o ,\c)$.
So, for each structure $\X$ chained by an ordinal $\b \in [\o , \c)$
by Theorem \ref{T012} there is an ordinal $\a$ such that $|\b|\leq \a \leq \b$ and $\B _\X \cong \B _\a$ and, hence
$$
\B _\X \cong \Col (\o _1 ,\c) \;\mbox{ or }\; \B _\X \cong \Col (\o ,\c).
$$
We will have $\B _\X \cong \Col (\o _1 ,\c)$
iff the ordinal $\a$ written in the form (\ref{EQ009}) satisfies the assumption of (i) of Theorem \ref{T226},
or $\X$ is a constant structure (note that, by Facts \ref{T4130} and \ref{T209}(a), $\ro (P(\o)/\Fin)\cong \Col (\o _1 ,\c)$).
\end{ex}
\begin{ex}\label{EX005}\rm
If $\X$ is a non-constant structure chained by an uncountable indecomposable ordinal $\o ^\d$,
then by Theorem \ref{T012}(c) we have $\B _\X \cong \B _{\o ^\d}$
and, clearly, $|X|=|\o ^\d|=|\d|\geq \o _1$. We regard Theorem \ref{T207} and recall that $\k :=\cf (\d)$.
\begin{itemize}\itemsep=-2mm
\item By Thm 5.2 of \cite{Kord1}, if (D) or (E) holds, $2^\k =2^{|\d|}$ and ($2^{<\k}=\k$  or $2^\k=\k ^+$),
      then $\B _{\o ^\d} \cong \Col (\o ,2^{|\d|} )$, and, hence, $\B _\X \cong \Col (\o ,2^{|X|} )$. \\
      For example, if  $2^{\o _1}=\o _2$ and $\X$ is a non-constant structure chained by the ordinal $\a =\o _2 ^{\o _2}$,
      then, since $\a =(\o ^{\o _2})^{\o _2}=\o ^{\o _2\o _2}$,
      we have $\d=\o_2\o_2$, so (E) holds and $|X|=|\d|=\cf (\d)=\k =\o _2=2^{<\o _2}$; thus
      \begin{equation}\label{EQ010}
      \B _{\o ^{\d}}=\B _{\o ^{\o_2\o_2}} \cong \B _\X \cong \Col (\o ,2^{|X|} )=\Col (\o ,2^{|\d|} )=\Col (\o ,2^{\o _2} ).
      \end{equation}
\item By Thm 5.7 of \cite{Kord1}, if $\P (\o ^\d)$ collapses $2^{|\d |}$ to $\o$,
      then $\B_{\o ^{\d +n}}\cong \Col (\o _1, 2^{|\d |})$, for all $n\in \N$. We note that here we have case (A) of Theorem \ref{T207}.\\
      For example, if  $2^{\o _1}=\o _2$ again and, for $n\in \N$, $\Y _n$ is a non-constant structure chained by the ordinal $\b _n =\o _2 ^{\o _2}\o^n$,
      then, since $\b _n=\o ^{\o _2\o _2 +n}=\o ^{\d +n}$, where $\d=\o_2\o_2$, by (\ref{EQ010}) for each $n\in \N$ we have
      \begin{equation}\label{EQ011}
      \B _{\o ^{\d +n}} \cong \B _{\Y _n} \cong \Col (\o_1 ,2^{|Y_n|} )=\Col (\o _1,2^{\o _2} ).
      \end{equation}
\end{itemize}
\end{ex}
\paragraph{Structures chained by uncountable real types}
The results of this paragraph confirm the important role of the Sacks perfect set forcing $\S$ in our context.
The following basic facts will be used in the sequel (see \cite{Kech}, pages 32, 83 and 89).

\begin{fac}\label{T013}
(a) If $A\in \Borel (\BR )\setminus [\BR]^{\leq \o}$, then $|A|=\c$ and $A$ contains a perfect set.

(b) If $f: \BR \rightarrow \BR$ is a Borel injection, then $f[\BR]\subset\BR$ is a Borel set.
\end{fac}
We remind the reader that if $\P$ and $\Q$ are partial orders,
then a function $e :\P \rightarrow \Q$ is called a {\it dense embedding},
we write $e :\P \hookrightarrow _d \Q$,
iff for each $p,q\in P$ we have:
{\rm (de1)} $p\leq q \Rightarrow e(p)\leq e(q)$,
{\rm (de2)} $p\perp q \Rightarrow e(p)\perp e(q)$,
{\rm (de3)} $e[P]$ is dense in $\Q$;
then we write $\P \hookrightarrow _d \Q$.
If $\tau$ is the topology on the set $P$ generated by the principal ideals $(\cdot ,p]$, $p\in P$,
the regular open algebra $\RO(\la P,\tau\ra)$ is denoted by $\ro (\P)$
and $e_\P : \P \hookrightarrow _d \ro (\P)$, where $e_\P (p)=(\cdot ,p]$, for all $p\in P$.
Moreover, $\ro (\P )$ is the  unique (up to isomorphism) complete Boolean algebra $\B$ such that $\P \hookrightarrow _d \B ^+$
(see \cite{Kun}, p.\ 63); the algebra $\ro (\P )$ is called the {\it Boolean completion of $\P$}.
We will use the following well known facts (see \cite{Kun}, pages  221 and 243).
\begin{fac}\label{T014}
$\P \mbox{ is a dense suborder of }\,\Q \;\Rightarrow\; \P \hookrightarrow _d \Q \;\Rightarrow\; \ro (\P) \cong \ro (\Q)$.
\end{fac}
\begin{te}\label{T004}
If $\,\BL$ is a suborder of the real line and contains a perfect set, then
$$
\B _\BL \cong \B _\BR \cong\ro (\S )\cong \ro(\Borel (\BR )/[\BR]^{\leq \o}).
$$
\end{te}
\dok
We divide the proof into a sequence of claims.\\[2mm]
\noindent
{\it Claim 1. If $f:\Q \hookrightarrow \BR$, then there is $F:\BR \hookrightarrow \BR$ such that $F\upharpoonright \Q =f$.}

\noindent
{\it Proof.}
For $x\in \BR$ there is $q_x\in \Q$ such that $x <q _x$, so if $\Q\ni q\leq x$ we have $f(q)<f(q _x)$ and, hence, there exists
\begin{equation}\label{EQ001}
F(x):=\sup \{ f(q):q\in \Q \land q\leq x\}.
\end{equation}
Thus if $p\in \Q$, then $F(p):=\sup \{ f(q):q\in \Q \land q\leq p\}=f(p)$ and $F\upharpoonright \Q =f$.
If $x<y$ and $p,q\in \Q$, where $x<p<q<y$,
Then $F(x)\leq f(p)<f(q)\leq F(y)$ and, hence, $F$ is strictly increasing.
\kdok
\noindent
{\it Claim 2. $\P (\BR)\subset \Borel (\BR)\setminus [\BR]^{\leq \o}$.}

\noindent
{\it Proof.}
If $f\in \Emb (\BR )$, $a \in R$, $f^{-1}[(-\infty ,a)]\neq \emptyset$ and $x<y\in f^{-1}[(-\infty ,a)]$,
then $f(x)<f(y)<a$ and, hence, $x\in f^{-1}[(-\infty ,a)]$.
Thus $f^{-1}[(-\infty ,a)]$ is an initial part of $\BR$,
that is a set of the form $(-\infty , b)$ or  $(-\infty , b]$ and, hence, a Borel set.
So $f$ is a Borel injection and by Fact \ref{T013}(b) the set $f[\BR ]$ is Borel.
\kdok
\noindent
{\it Claim 3. If $f:\Q \hookrightarrow \BR$, then $\BR \hookrightarrow \overline{f[\Q]}$.}

\noindent
{\it Proof.}
Let $f:\Q \hookrightarrow \BR$.
By Claim 1 for the function $F:\BR \rightarrow \BR$ given by (\ref{EQ001})
we have $F:\BR \hookrightarrow \BR$, $F\upharpoonright \Q =f$ and we prove that $F[\BR]\subset \overline{f[\Q]}$.
For $x\in \Q$ we have $F(x)=f(x)\in f[\Q]\subset \overline{f[\Q]}$.
If $x\in \BR \setminus \Q$, if $U$ is a neighborhood of $F(x)$ and $(F(x)-\frac{1}{n}, F(x)+\frac{1}{n})\subset U$,
then by (\ref{EQ001}) there is $q\in \Q$ such that $q<x$ and $F(x)-\frac{1}{n}<f(q)<F(x)$
and, hence, $U\cap f[\Q]\neq \emptyset$.
So $F(x)\in \overline{f[\Q]}$ again.
\kdok
\noindent
{\it Claim 4. For $A\subset \BR$ we have: $\BR \hookrightarrow A$ iff $A$ contains a perfect set.}

\noindent
{\it Proof.}
If $f:\BR \hookrightarrow A$, then by Claim 2 $f[\BR]\subset A$ is a Borel set of size $\c$,
and by Fact \ref{T013}(a) contains a perfect set.
Conversely, if $\S \ni P\subset A$, then by Fact \ref{T013}(a) $|P|=\c$.
Assuming that $\la P,<\ra$ is a scattered linear order
we would have $\o _1 \hookrightarrow P$ or $\o _1^* \hookrightarrow P$ (see \cite{Rosen}, p.\ 87),
which is impossible since $\o _1, \o _1^* \not\hookrightarrow \BR$.
Thus there is $f:\Q \hookrightarrow P$ and by Claim 3 $\BR \hookrightarrow \overline{f[\Q]}\subset \overline{P}=P$ (since $P$ is closed).
\kdok
\noindent
{\it Claim 5. $\B _\BR\cong \ro(\Borel (\BR)/[\BR]^{\leq \o})$.}

\noindent
{\it Proof.}
First we show that $\B_\BR \cong \ro(\la \Borel (\BR)\setminus [\BR]^{\leq \o}, \subset\ra)$.
By Fact \ref{T014} it is sufficient to prove that
$\P (\BR)$ is a dense subset of the poset $\la \Borel (\BR)\setminus [\BR]^{\leq \o}, \subset\ra$.
So, if $A\in \Borel (\BR)\setminus [\BR]^{\leq \o}$,
then $A$ contains a perfect set and by Claim 4 there is $f: \BR \hookrightarrow A$. Thus $f[\BR]\in \P (\BR)$ and $f[\BR]\subset A$, and we are done.

Second, $\ro(\la \Borel (\BR)\setminus [\BR]^{\leq \o}, \subset\ra)\cong \ro(\Borel (\BR)/[\BR]^{\leq \o})$.
Namely, using Fact \ref{T013} we check that $e:\Borel (\BR)\setminus [\BR]^{\leq \o}\hookrightarrow_d (\Borel (\BR)/[\BR]^{\leq \o})^+ $,
where $e(A)=[A] :=\{ B\in \Borel (\BR): |A \vartriangle B|\leq \o\}$ and then apply Fact \ref{T014}.
\kdok
\noindent
{\it Claim 6. $\ro(\S )\cong \ro(\Borel (\BR)/[\BR]^{\leq \o})$.}

\noindent
{\it Proof.}
By Claim 5 and Fact \ref{T014} it is sufficient to prove that
$\S$ is a dense subset of the poset $\la\Borel (\BR)\setminus [\BR]^{\leq \o}, \subset\ra$.
So, if $P\in \S$, then $P$ is a closed set of size $\c$
and, hence, $P\in \Borel (\BR)\setminus [\BR]^{\leq \o}$.
In addition, if $A\in \Borel (\BR)\setminus [\BR]^{\leq \o}$,
then there is $P\in \S$ such that $P\subset A$;
thus, $\S$ is dense in $\Borel (\BR)\setminus [\BR]^{\leq \o}$.
\kdok
\noindent
{\it Claim 7. $\B _\BL \cong \ro (\S )\cong \ro(\Borel (\BR )/[\BR]^{\leq \o})$.}

\noindent
{\it Proof.}
By Theorem 2.10 of \cite{KDif}, if $\X$ and $\Y$ are bi-embeddable structures of the same language
(i.e.\ $\X \rightleftarrows \Y$), then $\B_\X \cong \B _\Y$.
So, since $\BL \hookrightarrow \BR$ and, by Claim 4, $\BR \hookrightarrow \BL$,
we have $\B_\BL \cong \B _\BR$ and our claim follows from Claims 5 and 6.
\hfill $\Box$
\begin{te}\label{T005}
If $\X$ is a non-constant relational structure chainable by a real order type containing a perfect set, then
$$
\B _\X \cong \ro (\S )\cong \ro(\Borel (\BR )/[\BR]^{\leq \o}).
$$
\end{te}
\dok
W.l.o.g.\ we assume that $\BL \in \L_\X$ and that $\BL$ is a suborder of $\BR$ containing a perfect set.
By Theorems \ref{T018}(c) and \ref{T004} it is sufficient to prove that $\BL \in AI$.
So, if $\BL =\BI +\BF$, then by Claim 4 we have $\BR \hookrightarrow \BL$
and, hence, $\BR \hookrightarrow \BI$, which implies that $\BL \hookrightarrow \BI$,
or $\BR \hookrightarrow \BF$, which implies $\BL \hookrightarrow \BF$, and we are done.
\hfill $\Box$
\begin{rem}\label{R000}\rm
Concerning Theorems \ref{T004} and \ref{T005} we recall that

1. The class of suborders of $\BR$ containing a perfect set includes uncountable Borel and analytic sets
and, under the Axiom of projective determinacy, the class of uncountable projective sets (see \cite{Kech}, pages 83, 226 and 326).
Thus, in models of ZFC, in all of these situations, regarding Theorems \ref{T004} and \ref{T005} we will have $\B _\BL\cong\B _\X \cong \ro (\S)$.
(In fact, we can make $\BL$ adding anything to any perfect set and theorems will work.)
We note that, assuming the existence of an inaccessible cardinal, Solovay \cite{Sol} constructed a model of ZF $+\;\neg$AC
such that every set of reals of size $\c$ contains a perfect subset.

2. On the other hand, if $\BL \subset \BR $ is a Bernstein set (see \cite{Kech}, p.\ 48 for a construction in ZFC),
then neither $\BL$ nor its complement contains a perfect set.
Moreover, Dushnik and Miller \cite{Dus} constructed a dense suborder $\BL$ of $\BR$ of size $\c$
such that $\Emb (\BL)=\{ \id _L\}$ (see also \cite{Rosen}, p.\ 147).
Thus $\P (\BL)=\{ L\}$ and, hence, $\B _\BL \cong 2$.

3. Finally, we recall Remark \ref{R001}: if PFA holds and a structure $\X$ is chained by an $\aleph _1$-dense suborder of $\BR$,
then  we have $\B_\X \cong  \Col (\o ,\c)$.
\end{rem}
\paragraph{Acknowledgement.}
This research was supported by the Science Fund of the Republic of Serbia,
Program IDEAS, Grant No.\ 7750027:
{\it Set-theoretic, model-theoretic and Ramsey-theoretic
phenomena in mathematical structures: similarity and diversity}--SMART.


{\footnotesize
\begin{thebibliography}{99}
\bibitem{BPS}
      B.\ Balcar, J.\ Pelant, P.\ Simon,
      The space of ultrafilters on $\N$  covered by nowhere dense sets,
      Fund.\ Math.\ 110,1 (1980) 11--24.
\bibitem{Balc1}
      B.\ Balcar, P.\ Simon,
      Disjoint refinement,
      in: Handbook of Boolean algebras, Vol. 2, 333--388, North-Holland, Amsterdam, 1989.
\bibitem{Baum}
      J.\ E.\ Baumgartner,
      All $\aleph _1$-dense sets of reals can be isomorphic,
      Fund.\ Math.\ 79,2 (1973) 101-106.
\bibitem{Dus}
      B.\ Dushnik, E.\ W.\ Miller,
      Concerning similarity transformations of linearly ordered sets,
      Bull.\ Amer.\ Math.\ Soc.\ 46 (1940) 322--326.
\bibitem{Fra}
      R.\ Fra\"{\i}ss\'{e},
      Theory of relations, Revised edition, With an appendix by Norbert Sauer,
      Studies in Logic and the Foundations of Mathematics, 145.
      North-Holland, Amsterdam, (2000)
\bibitem{Fras}
      C.\ Frasnay,
      Quelques probl\'emes combinatoires concernant les ordres totaux et les relations monomorphes,
      Ann.\ Inst.\ Fourier (Grenoble) 15,2 (1965) 415--524.
\bibitem{Gib}
      P.\ C.\ Gibson, M.\ Pouzet, R.\ E.\ Woodrow,
      Relational structures having finitely many full-cardinality restrictions,
      Discrete Math.\ 291,1-3 (2005) 115--134.
\bibitem{Hodg1}
      W.\ Hodges, A.\ H.\ Lachlan, S.\ Shelah,
      Possible orderings of an indiscernible sequence,
      Bull.\ London Math.\ Soc.\ 9,2 (1977) 212--215.
\bibitem{Kech}
      A.\ S.\ Kechris,
      Classical descriptive set theory.
      Graduate Texts in Mathematics, 156. Springer-Verlag, New York, 1995.
\bibitem{Koj}
      M.\ Kojman, S.\ Shelah,
      Fallen cardinals,
      Dedicated to Petr Vopenka,
      Ann.\ Pure Appl.\ Logic 109,1-2 (2001) 117--129.
\bibitem{Kun}
      K.\ Kunen,
      Set theory. An introduction to independence proofs.
      Studies in Logic and the Foundations of Mathematics, 102. North-Holland Publishing Co., Amsterdam-New York, 1980.
\bibitem{Ktow}
      M.\ S.\ Kurili\'c,
      From $A_1$ to $D_5$: towards a forcing-related classification of relational structures,
      J.\ Symb.\ Log.\ 79,1 (2014) 279--295.
\bibitem{Kscatt}
      M.\ S.\ Kurili\'c,
      Posets of copies of countable scattered linear orders,
      Ann.\ Pure Appl.\ Logic 165,3 (2014) 895--912.
\bibitem{Kord}
      M.\ S.\ Kurili\'c,
      Forcing with copies of countable ordinals,
      Proc. Amer. Math. Soc. 143,4 (2015) 1771--1784.
\bibitem{KDif}
      M.\ S. Kurili\'c,
      Different similarities,
      Arch.\ Math.\ Logic 54,7--8 (2015) 839--859.
\bibitem{KVau}
      M.\ S. Kurili\'c,
      Vaught's conjecture for monomorphic theories,
      Ann.\ Pure Appl.\ Logic 170,8 (2019) 910--920.
\bibitem{Kord1}
      M.\ S.\ Kurili\'c,
      Forcing with copies of uncountable ordinals,
      https://arxiv.org/pdf/2401.00302
\bibitem{KCol}
      M.\ S.\ Kurili\'c,
      Reduced products of collapsing algebras,
      https://arxiv.org/pdf/2403.17930.pdf.
\bibitem{KurTod}
      M.\ S.\ Kurili\'c, S.\ Todor\v cevi\'c,
      Forcing by non-scattered sets,
      Ann.\ Pure Appl.\ Logic 163 (2012) 1299--1308.
\bibitem{Laf}
      C.\ Laflamme, M.\ Pouzet, R.\ Woodrow,
      Equimorphy: the case of chains,
      Arch.\ Math.\ Logic 56,7--8 (2017) 811--829.
\bibitem{Lav}
      R.\ Laver,
      On Fra\"{\i}ss\'{e}'s order type conjecture,
      Ann.\ of Math.\ 93,2 (1971) 89--111.
\bibitem{Moo}
      J.\ T.\ Moore,
      A five element basis for the uncountable linear orders,
      Ann.\ of Math.\ (2) 163,2 (2006) 669--688.
\bibitem{Pou}
      M.\ Pouzet,
      Application de la notion de relation presque-encha\^{\i}nable au d\'{e}nombrement des restrictions finies d'une relation,
      Z.\ Math.\ Logik Grundlagen Math.\ 27,4 (1981) 289--332.
\bibitem{Rosen}
      J.\ G.\ Rosenstein,
      Linear orderings,
      Pure and Applied Mathematics, 98, Academic Press, Inc.\ Harcourt Brace Jovanovich Publishers, New York-London, 1982.
\bibitem{She}
      S.\ Shelah,
      Power set modulo small, the singular of uncountable cofinality,
      J.\ Symbolic Logic 72,1 (2007) 226--242.
\bibitem{Sier}
      W.\ Sierpi\'nski,
      Cardinal and ordinal numbers. Second revised edition.
      Monografie Matematyczne, Vol.\ 34.\ Panstwowe Wydawnictwo Naukowe (PWN), Warsaw, 1965.
\bibitem{Sim}
      P.\ Simon,
      Sacks forcing collapses ${\mathfrak c}$ to ${\mathfrak b}$,
      Comment.\ Math.\ Univ.\ Carolin.\ 34,4 (1993) 707--710.
\bibitem{Sol}
      R.\ M.\ Solovay,
      A model of set-theory in which every set of reals is Lebesgue measurable,
      Ann.\ of Math.\ (2) 92 (1970) 1--56.
\end{thebibliography}
}
\end{document}